\documentclass[12pt,a4paper]{article}
\usepackage{textcomp}
\usepackage{listings}
\usepackage[space=true]{accsupp}
\newcommand{\copyablespace}{\BeginAccSupp{method=hex,unicode,ActualText=00A0}\EndAccSupp{}}
\usepackage[T1]{fontenc}
\usepackage[english]{babel}
\usepackage{hyperref}
\usepackage{graphicx}
\usepackage{amsmath}
\usepackage{amssymb}
\usepackage{mathtools}
\usepackage[table]{xcolor}
\usepackage{fancyhdr}
\usepackage{longtable}
\usepackage{multicol}
\usepackage{enumerate}
\usepackage{enumitem}
\setlist[itemize]{leftmargin=5.5mm}
\usepackage{multirow}
\usepackage{natbib}
\usepackage{color}
\definecolor{darkgreen}{rgb}{0.1, 0.5, 0.2}
\usepackage{caption}
\usepackage{subcaption}
\usepackage{xurl}
\usepackage[utf8]{inputenc}
\usepackage{tikz}
\usepackage{pgfplots}
\pgfplotsset{compat=newest}

\usepackage{framed}
\usepackage{booktabs}
\usepackage{caption}
\usepackage{float}
\usepackage{titlesec}
\usepackage{capt-of}
\usepackage{array}
\usepackage{arydshln}
\usepackage{appendix}
\newtheorem{theorem}{Theorem}
\newtheorem{lemma} {Lemma}

\newtheorem{example}{Example}

\newtheorem{condition}{DC}
\newtheorem{condition1}{MC}
\newtheorem{condition2}{SC}

\newenvironment{proof}[1][Proof]{\noindent\textbf{#1.} }{\ \rule{0.5em}{0.5em}}

\title{Comparative analysis of the existence and uniqueness conditions of parameter estimation in paired comparison models}
\author{\ L\'aszl\'o Gyarmati$^{1*}$, \'Eva Orb\'an-Mih\'alyk\'o$^{1}$, Csaba Mih\'alyk\'o$^{1}$}
\date{}
\begin{document}
\pagenumbering{arabic}

\maketitle
\begin{center}
$^{1}$ Department of Mathematics, University of Pannonia, 8200 Veszprém, Hungary\\
$^{*}$ Corresponding author, Department of Mathematics, University of Pannonia, 8200 Veszprém, Egyetem u. 10., Hungary; Email: gyarmati.laszlo@phd.mik.uni-pannon.hu, +3688624000/6109\\
\bigskip
Email: gyarmati.laszlo@phd.mik.uni-pannon.hu,
orban.eva@mik.uni-pannon.hu, mihalyko.csaba@mik.uni-pannon.hu \\

\end{center}
\begin{abstract}
In this paper paired comparison models with stochastic background are investigated. We focus on the models which allow three options for choice and the parameters are estimated by maximum likelihood method. The existence and uniqueness of the estimator is a key issue of the evaluation. In the case of two options, a necessary and sufficient condition is given by Ford in the Bradley-Terry model. We generalize this statement for the set of strictly log-concave distribution. Although in the case of three options necessary and sufficient condition is not known, there are two different sufficient conditions which are formulated in the literature. In this paper we generalize them, moreover we compare these conditions. Their capacities to indicate the existence of the maximum are analyzed by a large number of computer simulations. These simulations support that the new condition indicates the existence of the maximum much more frequently then the previously known ones.
\noindent

\end{abstract}


\renewcommand{\baselinestretch}{1.24} \normalsize

\section{Introduction}
Comparisons in pairs are frequently used in ranking and rating problems. They are mainly applied when scaling is very uncertain, but comparing the objects to the others can guarantee more definite results. The area of the possible applications is extremely large, some examples are the followings: education \citep{sahroni2016design, kosztyan2020analyzing},  sports \citep{cattelan2013dynamic,gyarmati2023aggregated,orban2022application},  information retrieval \citep{jeon2013revisiting, gyarmati2022incomplete}, energy supply \citep{trojanowski2021prospects},
financial sector \citep{montequin2020bradley}, management \citep{canco2021ahp}.
The most popular method is AHP (Analytic Hierarchy Process) elaborated by Saaty \citep{Saaty1977,saaty2004decision} and developed by others, see for example the detailed literature in \citep{liu2020review}. The method has lots of advantages: more than two options, several methods for evaluation, opportunity of incomplete comparisons, simple condition for the uniqueness of the evaluation \citep{BOZOKI}, possibility of multi-level decision \citep{rahman2021multi}, the concept of consistency \citep{brunelli2014introduction}. Nevertheless, due to the lack of stochastic background, the usual statistical tools, like confidence intervals, testing hypotheses are out of the possibilities. 

Fundamentally different models of paired comparisons are Thurstone motivated stochastic models. The basic concept is the idea of latent random variables, presented in  \citep{Thurstone1927}. Thurstone assumed Gauss distributed latent random variables and allowed two options in decisions, "worse" and "better". The method was modified: Gauss distribution was replaced by logistic distribution in \citep{BradleyTerry1952} and the model is called Bradley-Terry model (BTM). One of its main advantage are the simple mathematical formulae. Thurstone applied least squares method for parameter estimation, BTM applies maximum likelihood estimation and the not-complicated formulae allow quick numerical methods for solving optimization problems.
The existence and uniqueness of the optimizer is a key issue in the case of ML estimations; necessary and sufficient condition for it is proved in \citep{FORD}.
The model was generalized for three options ("worse", "equal" and "better") in \citep{glenn1960ties} for Gauss distribution and in \citep{rao1967ties} for logistic distribution. The latter paper applied maximum likelihood parameter estimation. Davidson made further modifications in the model concerning ties in \citep{davidson1970extending}. For more than 3 options we can find generalization in \citep{agresti1992analysis} in the case of Bradley-Terry model, and in \citep{OrbanMihalyko2019} in the case of Gauss distribution. In \citep{orban2019incomplete} it was proved that the models require the same conditions in order to be able to evaluate the data uniquely in the case of a broad set of cumulative distribution functions for the latent random variables: the strictly log-concave property of the probability density function is the crucial point of the uniqueness, while the assurance of the existence is hidden in the data structure. We mention that Gauss distribution and logistic distribution are included in the set of distributions having strictly log-concave probability density function. Note that, due to the probabilistic background, the Thurstone motivated models have the opportunity of building in the homefield or first-mover advantage \citep{hankin2020generalization}, testing hypotheses \citep{szabo2016study}, making forecasts \citep{mchale2011bradley}, therefore, they are worth investigating. 

In \cite{yan2016ranking}, the author analyzes the structure of the comparisons allowing both two and three options in choice. The author emphasizes that not only the structure of the graph made from the compared pairs but the results of the comparisons affect the existence of MLE. He makes some data perturbations in the cases where there are comparisons, but some results do not occur. By these perturbations, the zero data values become positive, and these positive value guarantee the strongly connected property of the directed graph constructed by the wins. But these perturbations modify the data structures, therefore, it would be better to avoid them.

In \citep{bong2022generalized}, the authors investigate BTM with two options and provide estimations for the probability of the existence of MLE. The authors turn to the condition of Ford to check whether MLE exists uniquely or not. As condition of Ford is necessary and sufficient condition, it indicates explicitly whether the MLE works or not. But in the case of other distributions and/or more than two options these investigations could not be performed due to the lack of necessary and sufficient condition for the existence and uniqueness of MLE. 

To continue their research, it would be conducive to have (necessary and) sufficient condition for the existence and uniqueness. To the best knowledge of the authors, there is no such theorem in the research literature, only two sufficient conditions is known. In this paper we compare the known conditions, we formulate their generalization, 
and we prove it. Then, we compare the applicability of the different conditions from the following point of view: how often and for what kind of parameters are they able to indicate the existence and uniqueness of MLE. We make large numbers of computer simulations and we use them to answer these questions. 

The paper is organised as follows: In Section \ref{IM} the investigated model is described. In Section  \ref{CEU} we present new conditions under which the existence and uniqueness is fulfilled. The proof can be found in Appendix A. In Section \ref{AC} the simulation results concerning the applicability are presented. Finally a short summary is given.

\section{The investigated model}\label{IM}
\begin{figure}[h]
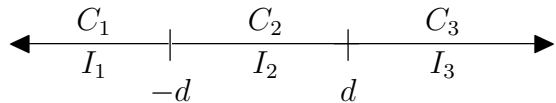
{}\[
{\tiny \blacktriangleleft}\hspace*{-0.5038pc}\dfrac{\hspace*{0.75cm}C_{1}%
\hspace*{0.75cm}}{I_{1}}\hspace*{-1.1082pc}%
\begin{array}
[c]{c}%
\\
\mid\\
-d%
\end{array}
\hspace*{-1.1082pc}\dfrac{\hspace*{1cm}C_{2}\hspace*{1cm}}{I_{2}}%
\hspace*{-1.1124pc}%
\begin{array}
[c]{c}%
\\
\mid\\
d%
\end{array}
\hspace*{-1.2112pc}\frac{\hspace*{1.2cm}C_{3}\hspace*{1.2cm}}{I_{3}}%
\hspace*{-1.2169pc}%
\hspace*{0.5pc}{\tiny \blacktriangleright}%
\
\]
\caption{The options and the intervals belonging to them }%
\label{INTFIG}%

\end{figure}

Let the number of the different objects to evaluate be denoted by $n$, and
let the objects be referred to as $1,2,...,n$.\ We want to evaluate them on the basis of the opinions of some persons called observers.
Let us denote the latent random variable belonging to the $i^{th}$ object by
$\xi_{i},$ $i=1,2,...,n$. Let the number of the options in a choice be $s=3$, namely "worse", "equal" and "better", denoted by $C_{1}$, $C_{2}$ and $C_{3}$. We split the set of the real
numbers $\mathbb{R}$ into $3$ \ intervals, which have no common elements. Each option in judgment corresponds to an interval in the
real line, the correspondence is noted by the same index. If the judgment
between the $i^{th}$ and $j^{th}$ objects is the option $C_{k},$ then we
assume that the difference $\xi_{i}-\xi_{j}$ of the latent random variables
$\xi_{i}$ and $\xi_{j}$ is in the interval $I_{k}, k=1,2,3.$
The intervals are determined by their initial points and endpoints, which are -$\infty$, $-d$, $d$ and $\infty$, $I_{1}$=(-$\infty$,$-d$), $I_{2}$=[$-d,d$] and  $I_{3}=$($d$,$\infty$).
The above intervals together with the corresponding options are presented in Figure \ref{INTFIG}.
We can write the differences of the latent random variables in the following form: 

\begin{equation}
\xi_{i}-\xi_{j}=m_{i}-m_{j}+\eta_{i,j}, i=1,...,n, j=1,...,n, i\neq j.
\label{kul}
\end{equation}
Now
\begin{equation}
E(\xi_{i})=m_{i}
\label{vhe}
\end{equation}
and $\eta _{i,j}$ are identically distributed random variables with expectation 0. The ranking of the expectations determines the ranking of the objects and the differences in their values give information concerning the differences of the strengths. We want to estimate the expectations and the value of the border of "equal" ($d$) on the basis of the data. For that we use maximum likelihood estimation.

The probabilities of the
events can be computed on the basis of the assumptions concerning the distributions of $\eta_{i,j}$  as follows:

\begin{equation}
P(\xi _{i}-\xi _{j}\in I_{1})=P(\xi _{i}-\xi _{j}<-d)=F(-d-(m_{i}-m_{j}))
\label{vveszit}
\end{equation}
\begin{equation}
P(\xi _{i}-\xi _{j}\in I_{2})=P(-d<=\xi _{i}-\xi _{j}<=d)=F(d-(m_{i}-m_{j}))-F(-d-(m_{i}-m_{j}))
\label{vdontetelen}
\end{equation}
\begin{equation}
P(\xi _{i}-\xi _{j}\in I_{3})=P(d<\xi _{i}-\xi _{j})=1-F(d-(m_{i}-m_{j}))
\label{vnyer}
\end{equation}
where $F$ is the (common) cumulative distribution function (c.d.f) of $\eta_{i,j}$.

Let the number of observers be $r$. The judgment produced by the $u^{th}$
observer ($u=1,2,...,r$) concerning the comparison of the $i^{th}$ and the
$j^{th}$ objects is encoded by the elements of a $4$ dimensional matrix which has
only 0 and 1 coordinates depending on the choice of the respondent. The third indices correspond to the options in choices, $k$=1,2,3 are for judgments "worse", "equal", and "better", respectively. The matrix of all judgments be $X,$ having 4 dimensions,
$i=1,2,...,n, j=1,2,...,n, k=$1, 2, 3, $u=1,2,...,r$ and

\rule[20pt]{0pt}{0pt}

$X_{i,j,k,u}$=$\left\{
\begin{array}
[c]{c}%
1,\text{if the opinion of the \ }u^{th\text{ }}\text{observer in pursuance
\ \ \ \ \ \ \ \ \ \ \ \ \ \ }\\
\text{of the comparison of the }i^{th}\text{ and the }j^{th}\text{ objects is
}C_{k}\\
\multicolumn{1}{l}{0,\text{ otherwise}}%
\end{array}
\right.$

\rule[20pt]{0pt}{0pt}

Let $X_{i,i,k,u}=0$. Of course, due to the symmetry, $X_{i,j,k,u}=X_{j,i,4-k,u}$. It expresses that if the $i^{th}$ object is "better" than the $j^{th}$ object, then the $j^{th}$ object is "worse" than the $i^{th}$ object, according to the judgment of the $u^{th} $ respondent.

Let $A_{i,j,k}=\sum_{u=1}^{r}X_{i,j,k,u}$ be the number of observations
$C_{k}$ in pursuance of the comparison of the $i^{th}$ and the $j^{th}$
objects and let $A$ denote the three dimensional matrix containing the
elements $A_{i,j,k}.$ Of course, $A_{i,j,k}=A_{j,i,4-k}.$

The likelihood function expresses the probability of the
sample in the function of the parameters. Assuming independent judgments, the likelihood function is%

\begin{equation}
L(X|m_{1},m_{2},...,m_{n},d)=%
{\textstyle\prod_{k=1}^{3}}
{\textstyle\prod_{i=1}^{n-1}}
{\textstyle\prod_{j=i+1}^{n}}
\left(  P(\xi _{i}-\xi _{j}\in I_{k})\right)  ^{A_{i,j,k}}%
\label{MLE}
\end{equation}
which has to be maximized in $\underline m=(m_{1},...,m_{n})$ and $0<d.$ 

One can realize that the likelihood function depends on the differences of the parameters $m_i$, therefore, one of them can be fixed.

\section{Conditions for the existence and uniqueness} \label{CEU}
In \citep{FORD}, the author presents a necessary and sufficient condition for the existence and uniqueness of MLE, if there are only two options for choice and $F$, the c.d.f. of $\eta_{i,j}$, is the logistic c.d.f.. The condition is the following: for arbitrary non-empty partition of the objects, $S$ and $\overline{S}$, there exists at least one element of $S$, which is "better" than an element of $\overline{S}$, and vice versa.
In \citep{davidson1970extending}, the author states that this condition supplemented with the condition "there is at least one tie ("equal")" is enough for having a unique maximizer in a modified Bradley-Terry model. The theorem assumes logistic distribution, its proof uses this special form, therefore, it is valid only for the investigated special model. Now we prove it for a broad set of c.d.f.'s. We require the following properties: $F$ is a c.d.f. with $0<F(x)<1$, $F$ is three times continuously differentiable, its probability density function $f$ is symmetric and the logarithm of $f$ is a strictly concave function in $\mathbb{R}$. Gauss and logistic distribution belong to this set, together with lots of others. Let us denote the set of these c.d.f.-s by $\mathbb{F}$. 

First we state the following generalization of Ford's theorem:
\begin{theorem}
Let $F\in \mathbb{F}$ and suppose that there are only two options in choice. Fix the value of the parameter $m_1=0$. The necessary and sufficient condition for the existence and uniqueness of MLE is the following: for arbitrary non-empty partition of the objects, $S$ and $\overline{S}$, there exists at least one element of $S$, which is "better" than an element of $\overline{S}$, and vice versa.
\label{GF}
\end{theorem}
The proof of sufficiency relies on the argumentation of Theorem \ref{TH3} omitting variable $d$. The used steps are (ST3), (ST5), and (ST6) in Appendix A. In the last step, the strictly concave property of $logL$ can be concluded from the theory of logarithmic concave measures \citep{prekopa1973logarithmic}. The necessity is obvious: if there would be a partition without "better" from one subset to another, then each element of this subset would be "worse" than the elements of the complement, but the measure of "worse" could not be estimated. The likelihood function would be monotone increasing, consequently, the maximum would not be reached. 

Returning to the case of three options, we formulate conditions of Davidson in the followings:
\begin{condition}
$\label{DC1}$
There exists an index pair $(i_{1},j_{1})$ for which $0<$A$%
_{i_{1},j_{1},{2}}.$
\end{condition}
\begin{condition}
For any non-empty partition of the objects $S$ and $\overline{S}$, there exists at least two index pairs ($i{_2}$,$j{_2}$) and ($i{_3}$,$j{_3}$)
$i_{2},i_{3}\in S$, $j_{2},j_{3}\in \overline{S}$
for which 
$0<A_{i_{2},j_{2},3}$ 
and 
$0<A_{i_3,j_3,1}.$ 
\label{DC2}
\end{condition}

Condition DC \ref{DC1} expresses that there is a judgment "equal". Condition DC \ref{DC2} coincides with the condition of Ford in \citep{FORD} in the case of two options. It expresses that there is at least one object in both subsets which is "better" than an object in the complement.

\begin{theorem}
Let $F\in \mathbb{F}$. If conditions DC \ref{DC1} and DC \ref{DC2} hold, then, fixing $m_{1}=0$, the likelihood function (\ref{MLE}) attains its maximal value and its argument is unique.
\label{THD}
\end{theorem}
Theorem \ref{THD} is the consequence of a more general statement, Theorem \ref{TH3}, which will be proved in Appendix \ref{A}. 

Now we turn to another set of conditions  which guarantees the existence and uniqueness of MLE. These conditions will be abbreviated by the initial letters MC.

\begin{condition1}
There is at least one index pair $(i_{1},j_{1})$ for which 

$0<A_{i_{1},j_{1},2}$ 
holds.
\label{MC1}
\end{condition1}

\begin{condition1}
There is at least one index pair $(i_{2},j_{2})$ for which 

$0<A_{i_{2},j_{2},1}$ and $0<A_{i_{2},j_{2},3}.$

\label{MC2}
\end{condition1}

Let us define the graph $G^{(M)}$ as follows: the nodes are the objects to be compared. There is an edge between two nodes $i$ and $j$, if 
$0<A_{i,j,2}$ or ($0<A_{i,j,1}$ and $0<A_{i,j,3}$)
hold.
\begin{condition1}
Graph $G^{(M)}$ is connected.
\label{MC3}
\end{condition1}

\begin{theorem}\citep{orban2019incomplete}
Let $F\in \mathbb{F}$. If conditions MC \ref {MC1}, MC \ref{MC2} and MC \ref{MC3} hold, then, after fixing $m_1$=0, the likelihood function (\ref{MLE}) attains its maximal value and the argument of the maximum is unique.
\label{THMIH}
\end{theorem}

\begin{figure}[H]
    \begin{minipage}[b]{0.45\textwidth}
        \centering
        \includegraphics[scale=0.25]{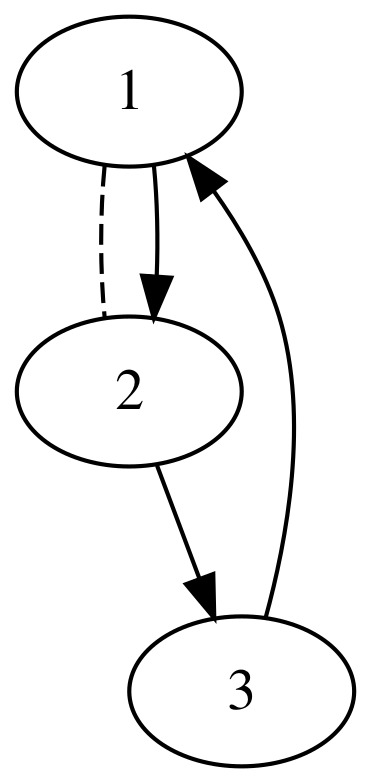}
                \caption{The graph of Example \Ref{EX1}.
                \newline
                - - - is for "equal",
                 -> is for "better".}
        \label{DIMN}
    \end{minipage}
    \hfill
    \begin{minipage}[b]{0.45\textwidth}
        \centering
        \includegraphics[scale=0.25]{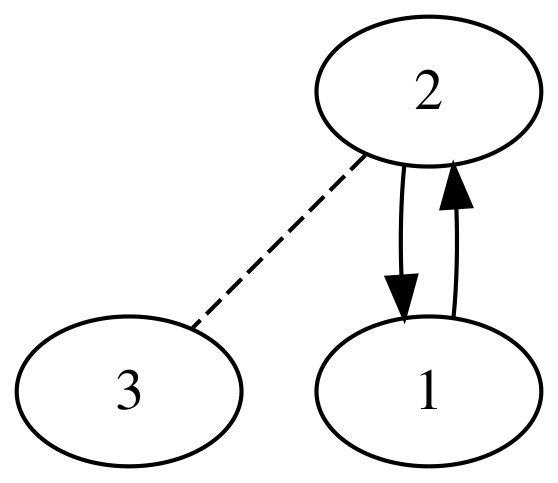}
         \caption{The graph of Example \Ref{EX2}.
         \newline   - - - is for "equal",
                 -> is for "better".}
         \label{DNMI}
    \end{minipage}
\end{figure}

To clear the relationship between conditions DC \ref{DC1}, DC \ref{DC2} and MC \ref{MC1}, MC \ref{MC2}, MC \ref{MC3} we present two examples.
In Example \ref{EX1}, DC \ref{DC1}, DC \ref{DC2} are satisfied but MC \ref{MC2} and MC \ref{MC3} are not. In Example \ref{EX2}, DC \ref{DC2} is not satisfied but MC \ref{MC1}, MC \ref{MC2}, MC  \ref{MC3} are. These examples expose that the sets of conditions DC and MC do not cover each other. Moreover, they support that MLE may exist uniquely even if DC \ref{DC1} and DC \ref{DC2} or MC \ref{MC1}, MC \ref{MC2} and MC \ref{MC3} do not hold. Therefore, we can see that neither conditions DC nor conditions MC are necessary conditions.
\begin{example}
Let n=3 and $A_{1,2,2}$=1, $A_{1,2,3}$=1, $A_{2,3,3}$=1, $A_{1,3,1}$=1 (see Figure \ref{DIMN}). Now both DC \ref{DC1} and DC \ref{DC2} hold, but MC \ref{MC3} does not.
\label{EX1}
\end{example}

\begin{example}
Let $n$=3 and $A_{1,2,1}$=1, $A_{1,2,3}$=1, $A_{2,3,2}$=1 (see Figure \ref{DNMI}). Now one can easily check that MC \ref{MC1}, MC \ref{MC2} and MC \ref{MC3} hold but DC \ref{DC2} does not.
\label{EX2}
\end{example}

The above theorems can be generalized. Let us introduce the following set of conditions denoted by SC:
\begin{condition2}
There is at least one index pair $(i_{1},j_{1})$ for which 

$0<A_{i_{1},j_{1},2}$ holds.
\label{SC1}
\end{condition2}
Let us introduce a graph belonging to the results of the comparisons as follows: let $DG^{(SC)}$ be a directed graph, the nodes are the objects, and there is a directed edge from $i$ to $j$ if there is an opinion according to which $i$ is "better" than $j$, that is $0<A_{i,j,3}$.
Now we can formulate the following conditions:
\begin{condition2}
There is a cycle in the directed graph $DG^{(SC)}.$
\label{SC2}
\end{condition2}

\begin{condition2}
For any non-empty partition of the objects $S$ and $\overline{S}$, there exists at least two (not necessarily different) index pairs ($i{_2}$,$j{_2}$) and ($i{_3}$,$j{_3}$)
$i_{2},i_{3}\in S$, $j_{2},j_{3}\in \overline{S}$  for which

\bigskip
$0<A_{i_{2},j_{2},3}$  
and $0<A_{i_{3},j_{3},1}$,
\bigskip

or there exists an index pair ($i{_4}$,$j{_4}$)
$i_{4}\in S$ and $j_{4}\in\overline{S}$
for which 
$0<A_{i_{4},j_{4},2}$.
\label{SC3}
\end{condition2}
It is easy to see that condition SC \ref{SC2} is more general than condition MC \ref{MC2} and condition SC \ref{SC3} is more general than condition DC \ref{DC2}. Condition SC \ref{SC3} expresses that any subset and its complement are interconnected by an opinion "better" or an opinion "equal". Here Condition  DC \ref{DC2} is replaced by a more general condition: next to "better" the opinion "equal" can also be appropriate judgment for connection.

To analyse the relationships between the sets of conditions DC, MC and SC we can recognize that 

(A) DC \ref{DC1}, MC \ref{MC1} and SC \ref{SC1} coincide. 

(B) If DC \ref {DC2} holds, then so does SC  \ref{SC2} and SC \ref{SC3}.

(C) If MC \ref{MC2} holds, so does SC \ref{SC2}.

(D) If MC \ref{MC3} holds, so does SC \ref{SC3}.

\begin{figure}[H]
    \begin{minipage}[b]{0.45\textwidth}
        \centering
        \includegraphics[scale=0.25]{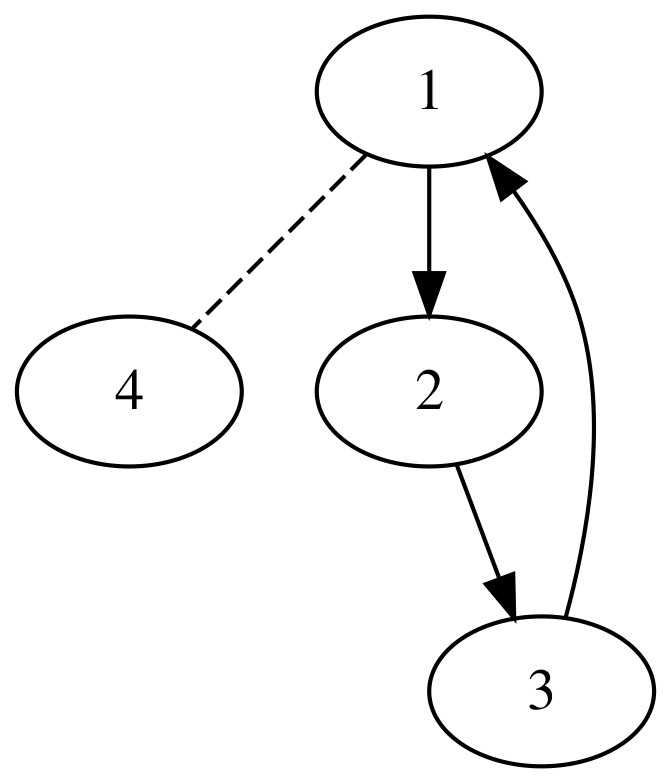}
        \caption{The graph of Example \Ref{EX3}.\newline
        - - - is for "equal",
                 -> is for "better".}
        \label{DNMN}
    \end{minipage}
    \hfill
    \begin{minipage}[b]{0.45\textwidth}
        \centering
        \includegraphics[scale=0.25
        ]{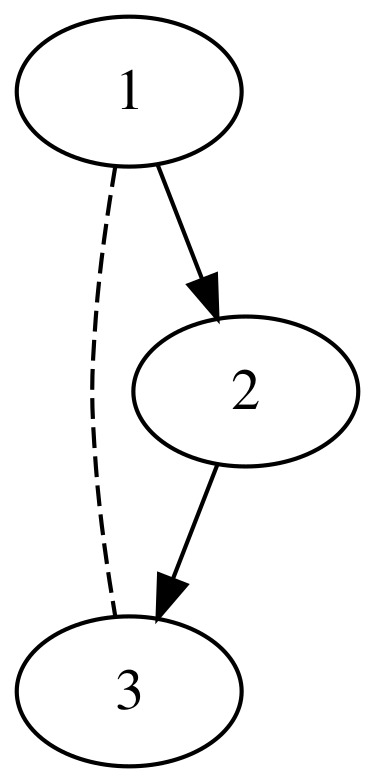}
         \caption{The graph of Example \Ref{EX4}.\newline
         - - - is for "equal",
                 -> is for "better".}
         \label{SNNSI}
    \end{minipage}
\end{figure}

These together present that conditions SC \ref{SC1}, SC \ref{SC2}, and SC \ref{SC3}
 are the generalization of the conditions DC and MC. 
 To show that SC is really a more general set of conditions we present Example \ref{EX3}.
 \begin{example}
 Let $n$=4, $A_{1,2,3}$=1, $A_{2,3,3}$=1, $A_{1,3,1}$=1 and $A_{1,4,2}$=1 (see Figure \ref{DNMN}). In this case neither condition DC \ref{DC2} nor MC \ref{MC2} hold, but SC \ref{SC1}, SC \ref{SC2} and SC \ref{SC3} do.  
 \label{EX3}
 \end{example}

Now we state the following theorem.
\begin{theorem}
Let $F\in \mathbb{F}$. If conditions SC \ref{SC1}, SC \ref{SC2} and SC \ref{SC3} hold, then, after fixing $m_1$=0, the likelihood function (\ref{MLE}) attains its maximum value and its argument is unique.

\label{TH3}
\end{theorem}
The proof of Theorem \ref{TH3} can be found in Appendix \ref{A}.

We note that Theorem \ref{THD} is a straightforward consequence of Theorem \ref{TH3}.

Unfortunately, conditions SC \ref{SC1}, SC \ref{SC2} and SC \ref{SC3} are not necessary conditions. One can prove that in the case of Example \ref{EX4} there exists a unique maximizer of function (\ref{MLE}) but SC \ref{SC2} does not hold.
\begin{example}
Let $n$=3, $A_{1,2,3}=1$, $A_{2,3,3}=1$ and $A_{1,3,2}=1$ (see Figure \ref{SNNSI}).
\label{EX4}
\end{example}

\section{Comparisons of the efficiency of the conditions} \label{AC}
In this section, we investigate in some special situations which sets of conditions  (conditions DC \ref{DC1}, DC \ref{DC2}; conditions MC \ref{MC1}, MC  \ref{MC2}, MC  \ref{MC3}; conditions SC \ref{SC1}, SC  \ref{MC2},  SC \ref{MC3}) are fulfilled, i.e. are able to detect the existence and the uniqueness of the maximizer. 

From the applications' perspective, there are such cases when the strengths of the objects to rank are close to each other and when they differ very much. On the other hand, there are such cases when the judgment "equal" is frequent, and such cases when it is rare. Referring to sports: in football and in chess the result draw comes up often, but in handball rarely. 

The most general set of conditions is the set SC. These conditions are fulfilled most frequently from the three sets of conditions. Nevertheless, it is interesting to what extent it is more applicable than the other two sets of conditions. For that we  made a large amount of computer simulations in the case of different parameter settings, and we investigated, how frequently the conditions are satisfied and how frequently we experience that the maximum exists.

We used Monte-Carlo simulation for the investigations. We fixed the differences between two expectations and the value of parameter $d$. This means that in our cases $\underline m=(0,h,2h,...,(n-1)h).$ We investigated 8 objects, and  we generated randomly  the pairs between which the comparisons exist. The number of comparisons was 8, 16, 32, 64. The results of the comparisons were also generated randomly, according to the probabilities  (\ref{vveszit}), (\ref{vdontetelen}) and (\ref{vnyer}). 

 In these random cases we checked whether conditions DC, MC, and SC are satisfied or not. Moreover we performed the numerical optimizations and we investigated whether the maximal value exists. We used 4 parameter ensembles, called situations, which are shown in Table \ref{tab:table_situations}.

In the presented situations, if the value of $h$ is small then the strengths of the objects are close to each other. It implies that many "better-worse" pairs could be formed during the simulations. On the other hand, if the value of $h$ is large, the strengths of the objects are far from each other, then we can expect only few "better-worse" pairs, but a great amount of "better" judgment. In terms of the number of "equal" judgments, if $d$ is large then lots of "equal" judgment could be formed during the simulations, while only few, when $d$ is small. The set of conditions DC can apply well the judgments "better", but it require only a single "equal" judgment. However, the set of conditions MC can use the judgments "equal" for connections, and the pairs "better-worse" judgments. Conditions SC do not require pairs, only judgments "better", in one circle. We recall that a single "better-worse" pair is appropriate as a circle. The judgments "equal" are well-applicable for this set of conditions, too. 

Table \ref{tab:table_situations} summarizes the situations with the presumable ratios of the "equal" judgments and "better-worse" pairs. In addition, Tables \ref{tab:situation_I}, \ref{tab:situation_II}, \ref{tab:situation_III} and \ref{tab:situation_IV} contains the numerical results of the simulations.
The order of the situations in terms of the number of the existence of the maximal values  is decreasing. Column MAX contains the number of the cases when the maximum exists. Columns DC/MAX, MC/MAX and SC/MAX present the ratios of the cases when the set of conditions DC, MC, SC hold, respectively. We can see that increasing the number of comparisons, the number of such cases when the maximal value exists and the ratios increase. We draw the attention to the fact that the values of the columns SC/MAX are less than 1 on several occasions. This detects again that SC is not a necessary condition.   

\begin{table}[h]\centering
\caption{Situations investigated} 
\begin{tabular}{|c|c|c|c|c|}
\hline
Situation&$h$&$d$&Rate of judgments "equal"& Rate of "better-worse" pairs \\\hline
I.&0.05&0.5&large&large\\\hline
II.&0.05&0.05&small&large\\\hline
III.&0.5&0.5&large&small\\\hline
IV.&0.5&0.05&small&small\\\hline

\end{tabular}
\label{tab:table_situations} 
\end{table}

We performed $10^8$ simulations per situation.
Table \ref{tab:situation_I} presents the results in Situation I. In this case we can see the DC/MAX rate is lower than the MC/MAX rate. We could predict it because there are lots of "equal" judgment. The SC/MAX rate is high even for 16 comparisons. In the case of 16 comparisons SC is 3.5 times better than MC and over 100 times better than DC. 

\begin{table}\centering
\caption{Situation I. ($h=0.05$, $d=0.5$)} 
\begin{tabular}{|c|c|c|c|c|}
\hline
Number of comparisons&MAX&DC/MAX&MC/MAX&SC/MAX\\\hline
8&57216&0&0.0921421&0.1941765\\\hline
16&38664325&0.0058568&0.2019802&0.7097257\\\hline

32&95920581&0.239853&0.8280385&0.9895364\\\hline

64&99987066&0.883599&0.9988596&0.9999986\\\hline

\end{tabular}
\label{tab:situation_I} 
\end{table}

Table \ref{tab:situation_II} presents the results of Situation II. In this case, the rate of "equal" is low, which does not favour the set of conditions MC. This is also reflected in the ratio MC/MAX, which is much worse than the ratio DC/MAX. The set of conditions SC still stands out among the other conditions.

\begin{table}\centering
\caption{Situation II. ($h=0.05, d=0.05$)} 
\begin{tabular}{|c|c|c|c|c|}
\hline
Number of comparisons&MAX&DC/MAX&MC/MAX&SC/MAX\\\hline
8&371&0&0&0.4070081\\\hline
16&5448890&0.3228876&0.0009119&0.9937707\\\hline

32&58963802&0.8708119&0.1898881&0.9999976\\\hline

64&92019027&0.9963352&0.9506307&1\\\hline

\end{tabular}
\label{tab:situation_II} 
\end{table}

Table \ref{tab:situation_III} shows the results of Situation III. Here the maximum values exist more rarely than in the previous two cases. In this case the number of "equal" decisions is high, while the number of "better-worse" pairs is low, which is favorable for the set of conditions MC and disadvantageous for the set of conditions DC, as we can see in Table \ref{tab:situation_III}. It can also be seen that none of the methods are as good as in the previous tables in terms of detecting the existence of the maximum. SC stands out again from the other two sets of conditions. Nevertheless, SC is able to show the existence of the maximum only in 73\% in the case of 32 comparisons, compared to 99\% in the previous situations. The set of conditions DC is almost useless, it is useful only in the cases 3.3\% even if the number of comparisons equals 64. The set of conditions MC method is slowly catching up and getting better, but for small numbers of comparisons (8, 16, 32) it is far from the much better SC criteria.

\begin{table}\centering
\caption{Situation III. ($h=0.5, d=0.5$)} 
\begin{tabular}{|c|c|c|c|c|}
\hline
Number of comparisons&MAX&DC/MAX&MC/MAX&SC/MAX\\\hline
8&248&0&0.0282258&0.0604839\\\hline
16&1025064&0.0005717&0.0532279&0.4203006\\\hline

32&23544050&0.004597&0.2771048&0.7256062\\\hline

64&76946023&0.0333163&0.8141669&0.95373\\\hline
\end{tabular}
\label{tab:situation_III} 
\end{table}

Table \ref{tab:situation_IV} presents the results in Situation IV. In the latter case, the numbers of "equal" choices and "better-worse" pairs are small, which is unfavorable MC, principally. In this situation, SC detects the existence of the maximal value exceptionally well. DC evinces them less fine, but it works better than MC. Nevertheless, for small numbers of comparisons, they are orders of magnitude weaker than SC.

\begin{table}\centering
\caption{Situation IV. ($h=0.5, d=0.05$)} 
\begin{tabular}{|c|c|c|c|c|}
\hline
Number of comparisons&MAX&DC/MAX&MC/MAX&SC/MAX\\\hline
8&2&0&0&1\\\hline
16&44246&0.1146209&0.00020355&0.9370956\\\hline
32&2621654&0.35173555&0.0184299&0.9965827\\\hline
64&25579173&0.6329823&0.37594615&0.99996685\\\hline
\end{tabular}
\label{tab:situation_IV} 
\end{table}

In all situations we have found that when we make few comparisons, SC is superior to the other conditions. As we make more and more comparisons, both other methods get better and better, but they are always worse than SC. The clear conclusion from the four tables is that the set of conditions SC is much more effective than the others, especially for small numbers of comparisons.

\section{Summary}
In this paper conditions guaranteeing the existence and uniqueness of the maximum likelihood parameter estimation are investigated. The case of general log-concave probability density function is studied. If two options are allowed, the usually applied Ford's condition is generalized from the logistic distribution to a wide set of distributions. This condition is necessary and sufficient condition. In the case of three options in decision, necessary and sufficient condition has not been proved, but there are two different sufficient conditions. We generalized them. A new set of conditions is proved which guarantees the existence and uniqueness of the maximizer. Moreover, we compare the conditions by the help of computer simulations and we have experienced that the set of the new conditions indicates the existence and uniqueness much more frequently, than the previously known conditions. Consequently, it provides more effective methods for such research which was preformed by Yan \citep{yan2016ranking} and Bong and Rinaldo \citep{bong2022generalized}.

The research includes the possibility of further developments. It would be desirable to set up the necessary and sufficient condition of the existence and uniqueness of the maximizer for the case of three options in choices, and simulations may help these findings. Further research is necessary to investigate the case of more than 3 options. These would be the subject of a next paper.

\bibliographystyle{apalike} 
\bibliography{main}

\begin{thebibliography}{}

\bibitem[Agresti, 1992]{agresti1992analysis}
Agresti, A. (1992).
\newblock Analysis of ordinal paired comparison data.
\newblock {\em Journal of the Royal Statistical Society: Series C (Applied
  Statistics)}, 41(2):287--297.

\bibitem[Bong and Rinaldo, 2022]{bong2022generalized}
Bong, H. and Rinaldo, A. (2022).
\newblock Generalized results for the existence and consistency of the mle in
  the {B}radley-{T}erry-{L}uce model.
\newblock In {\em International Conference on Machine Learning}, pages
  2160--2177. PMLR.

\bibitem[Boz{\'o}ki et~al., 2010]{BOZOKI}
Boz{\'o}ki, S., F{\"u}l{\"o}p, J., and R{\'o}nyai, L. (2010).
\newblock On optimal completion of incomplete pairwise comparison matrices.
\newblock {\em Mathematical and computer modelling}, 52(1-2):318--333.
\newblock \url{https://doi.org/10.1016/j.mcm.2010.02.047}.

\bibitem[Bradley and Terry, 1952]{BradleyTerry1952}
Bradley, R.~A. and Terry, M.~E. (1952).
\newblock Rank analysis of incomplete block designs: I. the method of paired
  comparisons.
\newblock {\em Biometrika}, 39(3/4):324--345.
\newblock \url{https://doi.org/10.2307/2334029}.

\bibitem[Brunelli, 2014]{brunelli2014introduction}
Brunelli, M. (2014).
\newblock {\em Introduction to the {A}nalytic {H}ierarchy {P}rocess}.
\newblock Springer.

\bibitem[Canco et~al., 2021]{canco2021ahp}
Canco, I., Kruja, D., and Iancu, T. (2021).
\newblock Ahp, a reliable method for quality decision making: A case study in
  business.
\newblock {\em Sustainability}, 13(24):13932.

\bibitem[Cattelan et~al., 2013]{cattelan2013dynamic}
Cattelan, M., Varin, C., and Firth, D. (2013).
\newblock Dynamic {B}radley--{T}erry modelling of sports tournaments.
\newblock {\em Journal of the Royal Statistical Society: Series C (Applied
  Statistics)}, 62(1):135--150.

\bibitem[Davidson, 1970]{davidson1970extending}
Davidson, R.~R. (1970).
\newblock On extending the {B}radley-{T}erry model to accommodate ties in
  paired comparison experiments.
\newblock {\em Journal of the American Statistical Association},
  65(329):317--328.

\bibitem[Ford~Jr, 1957]{FORD}
Ford~Jr, L.~R. (1957).
\newblock Solution of a ranking problem from binary comparisons.
\newblock {\em The American Mathematical Monthly}, 64(8P2):28--33.

\bibitem[Glenn and David, 1960]{glenn1960ties}
Glenn, W. and David, H. (1960).
\newblock Ties in paired-comparison experiments using a modified
  {T}hurstone-{M}osteller model.
\newblock {\em Biometrics}, 16(1):86--109.

\bibitem[Gyarmati et~al., 2022]{gyarmati2022incomplete}
Gyarmati, L., Orb{\'a}n-Mih{\'a}lyk{\'o}, {\'E}., Mih{\'a}lyk{\'o}, C.,
  Boz{\'o}ki, S., and Sz{\'a}doczki, Z. (2022).
\newblock The incomplete {A}nalytic {H}ierarchy {P}rocess and {B}radley-{T}erry
  model: (in)consistency and information retrieval.
\newblock {\em arXiv preprint arXiv:2210.03700}.

\bibitem[Gyarmati et~al., 2023]{gyarmati2023aggregated}
Gyarmati, L., Orb{\'a}n-Mih{\'a}lyk{\'o}, {\'E}., Mih{\'a}lyk{\'o}, C., and
  Vathy-Fogarassy, {\'A}. (2023).
\newblock Aggregated rankings of top leagues’ football teams: Application and
  comparison of different ranking methods.
\newblock {\em Applied Sciences}, 13(7):4556.

\bibitem[Hankin, 2020]{hankin2020generalization}
Hankin, R.~K. (2020).
\newblock A generalization of the {B}radley--{T}erry model for draws in chess
  with an application to collusion.
\newblock {\em Journal of Economic Behavior \& Organization}, 180:325--333.

\bibitem[Jeon and Kim, 2013]{jeon2013revisiting}
Jeon, J.-J. and Kim, Y. (2013).
\newblock Revisiting the {B}radley-{T}erry model and its application to
  information retrieval.
\newblock {\em Journal of the Korean Data and Information Science Society},
  24(5):1089--1099.

\bibitem[Koszty{\'a}n et~al., 2020]{kosztyan2020analyzing}
Koszty{\'a}n, Z.~T., Orb{\'a}n-Mih{\'a}lyk{\'o}, {\'E}., Mih{\'a}lyk{\'o}, C.,
  Cs{\'a}nyi, V.~V., and Telcs, A. (2020).
\newblock Analyzing and clustering students' application preferences in higher
  education.
\newblock {\em Journal of Applied Statistics}, 47(16):2961--2983.

\bibitem[Liu et~al., 2020]{liu2020review}
Liu, Y., Eckert, C.~M., and Earl, C. (2020).
\newblock A review of fuzzy ahp methods for decision-making with subjective
  judgements.
\newblock {\em Expert Systems with Applications}, 161:113738.

\bibitem[McHale and Morton, 2011]{mchale2011bradley}
McHale, I. and Morton, A. (2011).
\newblock A {B}radley-{T}erry type model for forecasting tennis match results.
\newblock {\em International Journal of Forecasting}, 27(2):619--630.

\bibitem[Montequ{\'\i}n et~al., 2020]{montequin2020bradley}
Montequ{\'\i}n, V.~R., Balsera, J. M.~V., Pilo{\~n}eta, M.~D., and P{\'e}rez,
  C.~{\'A}. (2020).
\newblock A {B}radley-{T}erry model-based approach to prioritize the balance
  scorecard driving factors: The case study of a financial software factory.
  mathematics, 8 (2).
\newblock {\em Applications of Operational Research and Mathematical Models in
  Management}, page 107.

\bibitem[Orb{\'a}n-Mih{\'a}lyk{\'o} et~al., 2022]{orban2022application}
Orb{\'a}n-Mih{\'a}lyk{\'o}, {\'E}., Mih{\'a}lyk{\'o}, C., and Gyarmati, L.
  (2022).
\newblock Application of the generalized {T}hurstone method for evaluations of
  sports tournaments’ results.
\newblock {\em Knowledge}, 2(1):157--166.

\bibitem[Orb\'an-Mih\'alyk\'o et~al., 2019a]{OrbanMihalyko2019}
Orb\'an-Mih\'alyk\'o, {\'E}., Mih\'alyk\'o, {\relax Cs}., and Koltay, L.
  (2019a).
\newblock A generalization of the {T}hurstone method for multiple choice and
  incomplete paired comparisons.
\newblock {\em Central European Journal of Operations Research},
  27(1):133--159.
\newblock \url{https://doi.org/10.1007/s10100-017-0495-6}.

\bibitem[Orb\'an-Mih\'alyk\'o et~al., 2019b]{orban2019incomplete}
Orb\'an-Mih\'alyk\'o, {\'E}., Mih\'alyk\'o, {\relax Cs}., and Koltay, L.
  (2019b).
\newblock Incomplete paired comparisons in case of multiple choice and general
  log-concave probability density functions.
\newblock {\em Central European Journal of Operations Research},
  27(2):515--532.

\bibitem[Pr{\'e}kopa, 1973]{prekopa1973logarithmic}
Pr{\'e}kopa, A. (1973).
\newblock On logarithmic concave measures and functions.
\newblock {\em Acta Scientiarum Mathematicarum}, 34:335--343.

\bibitem[Rahman et~al., 2021]{rahman2021multi}
Rahman, H.~U., Raza, M., Afsar, P., Alharbi, A., Ahmad, S., and Alyami, H.
  (2021).
\newblock Multi-criteria decision making model for application maintenance
  offshoring using analytic hierarchy process.
\newblock {\em Applied Sciences}, 11(18):8550.

\bibitem[Rao and Kupper, 1967]{rao1967ties}
Rao, P. and Kupper, L.~L. (1967).
\newblock Ties in paired-comparison experiments: A generalization of the
  {B}radley-{T}erry model.
\newblock {\em Journal of the American Statistical Association},
  62(317):194--204.

\bibitem[Saaty, 1977]{Saaty1977}
Saaty, T.~L. (1977).
\newblock A scaling method for priorities in hierarchical structures.
\newblock {\em Journal of Mathematical Psychology}, 15(3):234--281.
\newblock \url{https://doi.org/10.1016/0022-2496(77)90033-5}.

\bibitem[Saaty, 2004]{saaty2004decision}
Saaty, T.~L. (2004).
\newblock Decision making—the analytic hierarchy and network processes
  (ahp/anp).
\newblock {\em Journal of systems science and systems engineering},
  13(1):1--35.

\bibitem[Sahroni and Ariff, 2016]{sahroni2016design}
Sahroni, T.~R. and Ariff, H. (2016).
\newblock Design of analytical hierarchy process (ahp) for teaching and
  learning.
\newblock In {\em 2016 11th International Conference on Knowledge, Information
  and Creativity Support Systems (KICSS)}, pages 1--4. IEEE.

\bibitem[Szab{\'o} et~al., 2016]{szabo2016study}
Szab{\'o}, F., K{\'e}ri, R., Schanda, J., Csuti, P., and
  Mih{\'a}lyk{\'o}-Orb{\'a}n, E. (2016).
\newblock A study of preferred colour rendering of light sources: Home
  lighting.
\newblock {\em Lighting Research \& Technology}, 48(2):103--125.

\bibitem[Thurstone, 1927]{Thurstone1927}
Thurstone, L. (1927).
\newblock A law of comparative judgment.
\newblock {\em Psychological Review}, 34(4):273–286.
\newblock \url{https://doi.org/10.1037/h0070288}.

\bibitem[Trojanowski and Kazibudzki, 2021]{trojanowski2021prospects}
Trojanowski, T.~W. and Kazibudzki, P.~T. (2021).
\newblock Prospects and constraints of sustainable marketing mix development
  for {P}oland’s high-energy consumer goods.
\newblock {\em Energies}, 14(24):8437.

\bibitem[Yan, 2016]{yan2016ranking}
Yan, T. (2016).
\newblock Ranking in the generalized {B}radley--{T}erry models when the strong
  connection condition fails.
\newblock {\em Communications in Statistics-Theory and Methods},
  45(2):340--353.

\end{thebibliography}
\addcontentsline{toc}{section}{References}
\bigskip
\bigskip
\bigskip
\LARGE
\appendix{Appendix} \label{APP}
\normalsize
\section{Proof of Theorem \ref{TH3}} \label{A}
\label{A1}

\begin{proof} First we mention that instead of (\ref{MLE}), its logarithm, the log-likelihood function
\begin{equation}
log L(X|m_{1},m_{2},...,m_{n},d)=%
{\textstyle\sum_{k=1}^{3}}
{\textstyle\sum_{i=1}^{n-1}}
{\textstyle\sum_{j=i+1}^{n}}
{A_{i,j,k}}\cdot log P(\xi_{i}- \xi_{j} \in I_{k}))  %
\label{LOGMLE1}
\end{equation} 
that is 
\begin{equation}
log L(X|m_{1},m_{2},...,m_{n},d)=%
0.5 \cdot {\textstyle\sum_{k=1}^{3}}
{\textstyle\sum_{i=1}^{n}}
{\textstyle\sum_{j=1}^{n}}
{A_{i,j,k}}\cdot log P(\xi_{i}- \xi_{j} \in I_{k}))  %
\label{LOGMLE}
\end{equation}
is maximized under the conditions 0<$d$ and $m_1$=0. We prove that (\ref{LOGMLE}) attains its maximal value under the conditions 0<$d$ and $m_1$=0 and the argument of the maximal value is unique. 

The steps of the proof are denoted by (ST1), (ST2), (ST3), (ST4), (ST5) and (ST6).

Computing the value of the log-likelihood function at $m=(0,0,0,...,0), d=1$ and denoting this by $logL_0$, the maximum has to be seek in such regions where the values of (\ref{LOGMLE}) are at least $logL_0$. Moreover, we note that every term of the sum in (\ref{LOGMLE}) is negative (or zero if $A_{i,j,k}=0$), consequently, the maximum can not be attained in those regions where any term is under $logL_0$. By investigating the limits of the terms, we will check which parameters can be restricted into a closed bounded regions. The proof of the existence relies on the Weierstass theorem: we restrict the range of $d$ and $m_2$,...,$m_n$ to bounded closed sets where the continuous function (\ref{LOGMLE}) has maximal value. For that, we proof some lemmas.

(ST1) The first step is to find a positive lower bound for the variable $d$.
\begin{lemma} 
  \label{L1}
  Condition SC \ref{SC1} guarantees that the maximum can be attained if $\varepsilon\leq d$ with an appropriate value of $0<\varepsilon$.
  \end{lemma} 
\begin{proof}
SC \ref{SC1} guarantees that there exists an index pair $i,j$ for which $0<A_{i,j,2}$. Now,
\begin{equation}
A_{i,j,2} \cdot log(F(d-(m_{i}-m_{j}))-F(-d-(m_{i}-m_{j})))\longrightarrow-\infty \quad if  \quad d\longrightarrow0. 
\end{equation}
If  $d\longrightarrow0$, the arguments of the c.d.f. tend to the same value, their difference tends to zero. Consequently, its logarithm, with a positive multiplier, tends to minus infinity.

As $0.5 \cdot A_{i,j,2} \cdot log(F(d-(m_{i}-m_{j}))-F(-d-(m_{i}-m_{j})))$<$logL_0$, if $d<\varepsilon$, we can restrict the region of $d$ to the subset $\varepsilon\leq d$ with an appropriate value of $0<\varepsilon$, while seeking the maximum.
\end{proof}

\rule[30pt]{0pt}{0pt}
(ST2) The next step is to find an upper bound for the variable $d$.
\begin{lemma}
\label{L2}
If $0<A_{i,j,3}$ then there exists an upper bound $K_{i,j}$ for which it holds that the maximum can be attained in the region  $d-(m_i-m_j)\leq K_{i,j}$.
\end{lemma}

\begin{proof}
It is easy to see that if $0<A_{i,j,3}$, then
\begin{equation}
\label{korlgy}
   A_{i,j,3} \cdot log(1-F(d-(m_{i}-m_{j})))\longrightarrow -\infty \textit{ supposing }  d-(m_{i}-m_{j})\longrightarrow \infty.
\end{equation} 
Consequently, there exists a value $K_{i,j}$ with the following property: if  $K_{i,j}<d-(m_{i}-m_{j})$ then $0.5 \cdot A_{i,j,3} \cdot log(1-F(d-(m_{i}-m_{j}))) < logL_0$, the maximum has to be seek on the region $d-(m_{i}-m_{j}) \leq K_{i,j}.$ It means that the maximum can be reached only in such regions where $d-(m_{i}-m_{j})$ has an upper bound.
\end{proof}
 
 Condition SC \ref{SC2} guarantees that there is a cycle ($i_1,i_2,...,i_h,i_1$) with directed edges from $i_k$ to $i_{k+1}$ $k=1,2,...,h$ and from $i_h$ to $i_1$ and these directed edges arise from $0<A_{i_k,i_{k+1},3}$ $k=1,...,h$ and $0<A_{i_h,i_1,3}$. We can assume that $i_1$=1. 
Lemma \ref{L2} implies that 
\begin{equation}
\label{ineq}
d-(m_{i_k}-m_{i_{k+1}})\leq K_{k,k+1} 
\textit{ and } 
d-(m_{i_h}-m_{i_{1}})\leq K_{h,1}.
\end{equation}
Using a common upper bound $K$ ($K_{k,k+1}\leq K$ and $K_{h,1}\leq K$), moreover, summing the inequalities in (\ref{ineq}), we get that $h\cdot d\leq\sum_{k=1}^{h} K.$ This proves that it is enough to seek the maximum in a closed bounded set of $d$.

\rule[30pt]{0pt}{0pt}
(ST3) Now let us turn to the upper and lower bounds of the parameters $m_i$.

Let us defined a graph $G^{(SC)}$ as follows: the vertices are the objects. There is a directed edge from $i$ to $j$ if $0<A_{i,j,3}$ ($i$ is "better" than $j$ according at least one opinion). There is a directed edge from $i$ to $j$ and also from $j$ to $i$ if $0<A_{i,j,2}$ (they are "equal" according at least one opinion).
We will use the following well-known statement. Condition SC \ref{SC3} is equivalent to the following condition: between any pair of objects $i$ and $j$ there is a directed path in $G^{(SC)}$ from one to the other.
\begin{lemma}
\label{L3}
If $0<A_{i,j,3}$ and $m_i \leq K_i$, then there exists an upper bound of $m_j$ denoted by $K_j$, with the following property: $0.5 \cdot A_{i,j,3} \cdot  log(1-F(d-(m_{i}-m_{j})))<logL_0$ if $K_j<m_j$, that is the maximum can be attained if $m_j \leq K_j$.
\end{lemma}
\begin{proof}
Recalling (\ref{korlgy}), we can conclude that $d-(m_i-m_j) \leq K_{i,j}.$ As $m_i \leq K_i$ and $0<d<K^{(d)}$, we get that $m_j \leq K_{i,j}+K_i.$
\end{proof}

We can interpret Lemma \ref{L3}, that the property "having upper bound" spreads in the direction of the edge "better" defined by $0<A_{i,j,3}$.

\begin{lemma}
\label{L4}
If $0<A_{i,j,1}$ (there is at least one opinion according to $i$ is "worse" than $j$) and the inequality $-B_i \leq m_i$ holds, then there exists a lower bound of $m_j$, denoted by $-B_j$, with the following property: $0.5 \cdot A_{i,j,1} \cdot logF(-d-(m_{i}-m_{j}))<logL_0$ if $m_j<-B_j$, that is the maximum can be attained if $-B_j \leq m_j$.
\end{lemma}
\begin{proof}
The statement is the straightforward consequence of the following: if $0<A_{i,j,1}$, then
\begin{equation}
    logF(-d-(m_{i}-m_{j}))\longrightarrow -\infty \textit{ supposing }  -d-(m_{i}-m_{j})\longrightarrow -\infty.
    \end{equation}

    \end{proof}

We can interpret Lemma \ref{L4}, that the property "having lower bound" spreads along the opinion "worse".

\rule[10pt]{0pt}{0pt}

(ST4) Finally, we investigate the effect of the existence of a opinion "equal" for the property boundedness.

\begin{lemma}
    \label{L5}
    Suppose that the parameter $d$ is bounded. If 0$<A_{i,j,2}$ and $m_i \leq U_i$, there exists an upper bound $U_j$ for which if $U_j<m_j$ then (\ref{LOGMLE}) $< logL_0$. It means, that the maximum has to be sought in the region $m_j \leq U_j$.

    If 0$<A_{i,j,2}$ and $-H_i \leq m_i$, there exists a lower bound $-H_j$ for which if $m_j<-H_j$ then (\ref{LOGMLE}) $< logL_0$. It means, that the maximum has to be sought in the region $-H_j\leq m_j$.
    \end{lemma}
    \begin{proof}
    It is easy to see that 
    \begin{equation}
        \lim_{-d-(m_{i}-m_{j})\rightarrow -\infty }A_{i,j,2} \cdot log(F(d-(m_{i}-m_{j}))-F(-d-(m_{i}-m_{j})))=-\infty,
    \label{korl1}
\end{equation}
and
\begin{equation}
     \lim_{d-(m_{i}-m_{j})\rightarrow \infty }A_{i,j,2} \cdot log(F(d-(m_{i}-m_{j}))-F(-d-(m_{i}-m_{j})))=-\infty.
    \label{korl2}
\end{equation}
Consequently, the maximum has to be in the following region:
\begin{equation}
-H_{i,j} \leq -d-(m_i-m_j)  
\textit{ and }
d-(m_i-m_j) \leq B_{i,j},
\end{equation}
respectively, with an appropriate bound $-H_{i,j}$ and $B_{i,j}$.
As $\epsilon \leq d \leq K^{(d)}, m_i \leq B_i$ implies
$m_j \leq$ $B_{i,j}+B_i$ and $ -H_{i,j}-H_i \leq m_j$. 
\end{proof}

We can summarize Lemma \ref{L5}, that a both properties "having upper bound" and "having lower bound" spread with opinion "equal". It behaves as a "better" and a "worse" opinion, altogether. 

\rule[10pt]{0pt}{0pt}

(ST5)
Now we can prove that it is enough to seek the maximum on closed bounded set of every parameter $m_i$. Starting out of $m_1=0$, there exists a directed path from {1} to {$i$} in $G^{(SC)}$, along the edges defined by $0<A_{i,j,3}$ and  $0<A_{i,j,2}$. Walking along this path, and recalling that $m_1=0$, the property "having upper bound" spreads from $1$ to object $i$. The directed path from the object $i$ to $1$ is a reverse directed path from $1$ to $i$, and the property "having lower bound" of the object $1$ spreads to $i$ for every index $i$, consequently the expectations can be restricted into a closed bounded set. The maximal value of (\ref{LOGMLE}) can only be in these regions. As (\ref{LOGMLE}) is a continuous function, the Weierstrass theorem implies the existence of the maximal value.

\rule[30pt]{0pt}{0pt}
(ST6) The uniqueness of the argument at the maximal value is the consequence of the strictly concave property of the logarithm of the p.d.f.. Lemma 6  in \cite{orban2019incomplete}) implies the strictly concave property of the function (\ref{LOGMLE}) in $d-(m_i-m_j)$ and $-d-(m_i-m_j)$ for every index pair $(i,j)$ for which $0<A_{i,j,2}$ and in $d-(m_i-m_j)$ if $0<A_{i,j,3}$. 
Walking along the circle $DG^{(SC)}$ defined by SC \ref{SC2} and summing the arguments, we get that the function (\ref{LOGMLE}) is the strictly concave function of the parameter $d$. 

Now let us turn to the strictly concave property of (\ref{LOGMLE}) in the parameters $m_i$, $i=2,3,...,n$. There is a directed path from $1$ to $i$ defined by $0<A_{i,j,2}$ and $0<A_{i,j,3}$ in graph $G^{(SC)}$. Walking along it, we can conclude the strictly concave property of (\ref{LOGMLE}) in $d-(m_{i_k}-m_{i_{k+1}})$. Summing the arguments of the terms belonging to the path, we get, that the log-likelihood function is strictly concave in $ l\cdot d +m_i$, where $0<l$ is the length of the path. This fact and the strictly concave property in $d$ guarantee the strictly concave property in $m_{i}$. We get that function (\ref{LOGMLE}) is strictly concave in its every variable $m_i$ and $d$, therefore, the argument of the maximum has to be unique. 
\end{proof}

\end{document}